\documentclass{ifacconf}

\usepackage{graphicx}      
\usepackage{natbib}        

\begin{document}
\pagestyle{plain}
\pagenumbering{arabic}

\begin{frontmatter}

\title{On Some Geometric Interpretations \\ of Fractional-Order Operators} 

\thanks[footnoteinfo]{
This work is supported in parts by grants APVV-22-0508, VEGA 1/0674/23, 
ARO W911NF-22-1-0264. 
}

\author{Igor Podlubny} 

\address{BERG Faculty, Technical University of Kosice, 
   04200 Kosice, Slovakia  (e-mail: igor.podlubny@tuke.sk).}

\begin{abstract}                
This discussion paper presents some parts of the work in progress. 
It is shown that 
G.W. Leibniz was the first who raised the question about geometric interpretation
of fractional-order operators. 
Geometric interpretations of the Riemann--Liouville fractional integral 
and the Stieltjes integral are explained. 
Then, for the first time, a geometric interpretation of the Stieltjes derivatives
is introduced, which holds also for so-called ``fractal derivatives'', 
which are a particular case of Stieltjes derivatives. 

\end{abstract}


\end{frontmatter}

\section{Introduction: \\ "What is it in Geometry?"}

Based on available facts and the study of  G.W.~Leibniz's published correspondence, 
it is possible to state with confidence that he was the first who  opened the question about 
geometric interpretations of fractional-order differentiation and fractional-order 
integration. 

There are at least two available versions of the Leibniz's reply to G. F. de L'Hospital. 
The letter dated  September 30th, 1695,  first appeared in Leibniz's collected works published in 1849 (see \cite{Leibniz-1849}), and and many years after that in the more complete collection in 2004 (see \cite{Leibniz-2004}). 

This letter is frequently quoted because of the following Leibniz's prediction of the future importance of fractional differentiation (as first appeared in printed form in~\hbox{\cite{Leibniz-1849}}):

 \vspace*{-2ex}

\begin{center}
\includegraphics[width=8.4cm]{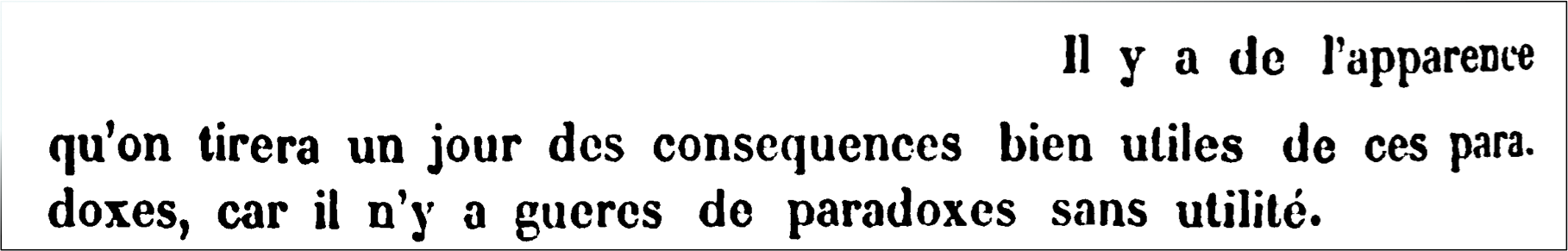}    
\end{center}
\begin{quote}
\emph{``It looks like one day we will draw very useful consequences from these paradoxes, because there are hardly any paradoxes without utility.''}
\end{quote}

\vspace*{1ex}

\begin{figure}[h!]
\begin{center}
\includegraphics[width=8.4cm]{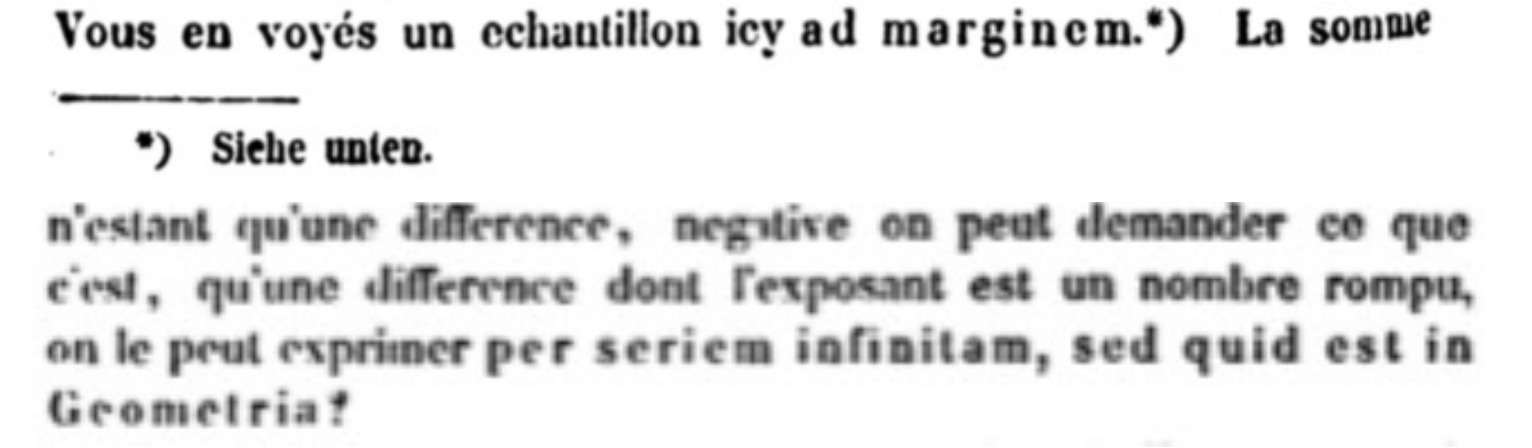}    
\includegraphics[width=8.4cm]{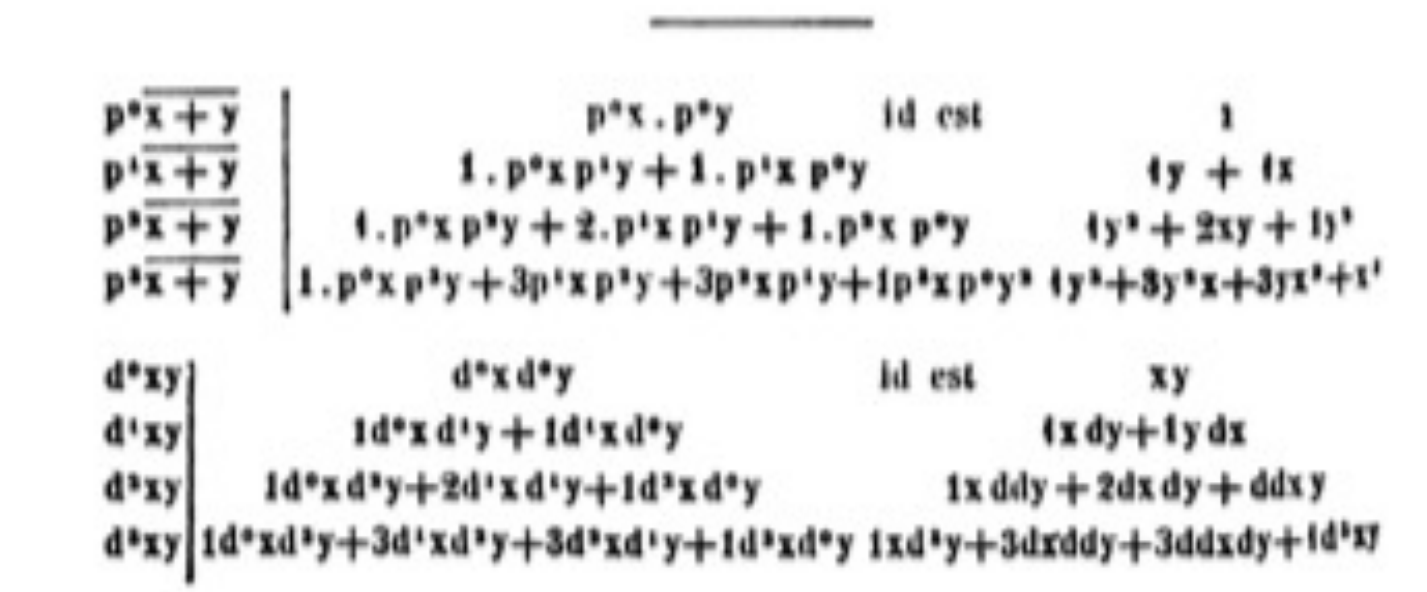} 
\caption{From Leibniz's letter as published  in \cite{Leibniz-1849}.} 
\label{fig:Leibniz-1849-quotation}
\end{center}
\end{figure}

In the shadow of this predictive statement, another important Leibniz's remark did not get proper attention. Namely, in the upper part of Fig.~\ref{fig:Leibniz-1849-quotation} the text of his letter is shown, referring to the remark on the margin (shown in the lower part):
 \begin{quotation}
\emph{``The sum is just a differential with negative exponent, and we can ask what is a differential whose exponent is a fractional number; we can express it by infinite series,  but what is it in Geometry?''}
\end{quotation}

\vbox{
This is where Leibniz raised the long-standing question: fractional derivatives and integrals represented by infinite series (or other formulas) can be used, but can they have a geometric interpretation? }

\begin{figure}[h!]
\begin{center}
\includegraphics[width=8.4cm]{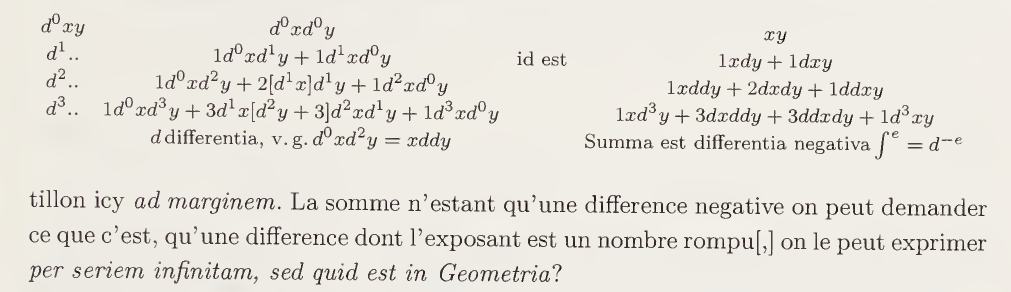}    
\caption{From Leibniz's letter as published  in \cite{Leibniz-2004}.} 
\label{fig:Leibniz-2004-quotation}
\end{center}
\end{figure}

In \cite{Leibniz-2004}, however, more details from the handwritten letter were added -- see Fig.~\ref{fig:Leibniz-2004-quotation}, and the sentence 
 \begin{quotation}
\emph{``Summa est differentia negativa $\int^e = d^{-e}$''} 
 \end{quotation}
along with the used notation and his remark on expressing fractional-order derivatives 
as series indicates that Leibniz realized that derivatives and integrals of any order 
can be considered as the same operation.

\section{Geometric interpretation of the Riemann--Liouville fractional integral}
It took long time -- more than 300 years -- until a widely accepted geometric interpretation 
of fractional-order integration was developed by \cite{Podlubny-2002}. 

\begin{figure}[h!]
\begin{center}
\includegraphics[width=8.4cm]{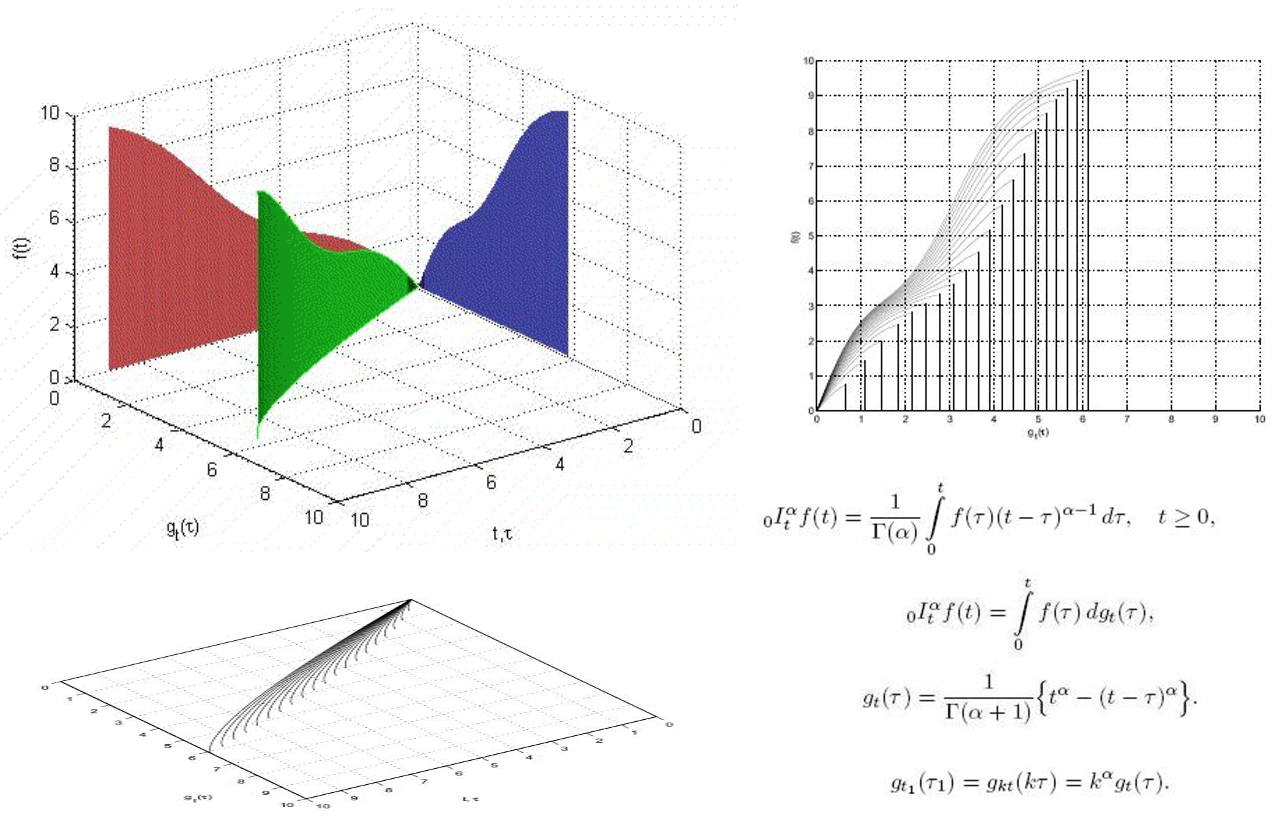}    
\caption{Geometric interpretation by \cite{Podlubny-2002}.} 
\label{fig:Podlubny-interpretation}
\end{center}
\end{figure}

In \cite{Podlubny-2002}, in order to utilize the idea by \cite{Bullock},
the Riemann--Liouville fractional integral is written 
in the form similar to the Stieltjes integral. 
This allowed to interpret the fractional-order integral 
as a ``live shadow'' (blue shape in Fig.~\ref{fig:Podlubny-interpretation}) 
of the ``live fence'' (green shape) of the height $f(t)$
along the curve 
$g_t(\tau) = \left\{ t^\alpha - (t-\tau)^\alpha    \right\} / \Gamma(\alpha+1)$, 
which is changing in a self-similar (self-scaling) way. 
The red area in Fig.~\ref{fig:Podlubny-interpretation} 
is just the classical geometric interpretation of the integral of the function $f(t)$, 
$\int_a^b f(t) dt$.

\section{Geometric interpretation of Stieltjes integral and Stieltjes derivative}

Looking back at the Stieltjes integral, we still can find more inspiration 
for interpretation of operators involving fractional orders.

\subsection{Stieltjes integral}

\vspace*{-1ex}

The integral of a function $f(x)$ with respect to another function $g(x)$,
\vspace*{-2ex}
\begin{equation}\label{Stieltjes-integral}
	\int_{a}^{b}
	f(t) \, dg(t),
\end{equation}

\vspace*{-1ex}
\noindent
was introduced by \cite{Stieltjes} as an auxiliary tool in his extensive study on continued fractions. 
Although the Stieltjes integral quickly found numerous applications in many fields of mathematics and physics, 
it took 94~years (almost a century!) until \cite{Bullock} suggested its geometric interpretation, 
which can be summarized as follows (see Fig.~\ref{fig:Stieltjes-interpretation}). 

Introduce an additional axis $\tau$ and build the ``fence'' of the height $y = f(t)$ 
along the curve $\tau = g(t)$, $a \leq t \leq b$ (green fence in Fig.~\ref{fig:Stieltjes-interpretation}). 
Projection of this fence on the plane $(t,y)$, shown in red color, is the classical geometric interpretation 
of the integral of the function $f(t)$, $\int_a^b f(t) dt$, 
while projection of this fence on the plane $(\tau,y)$, shown in blue color, 
gives  the geometric interpretation of the classical Stieltjes integral (\ref{Stieltjes-integral}).

\subsection{Stieltjes derivative}

\vspace*{-1ex}

Much less known is the notion of the Stieltjes derivative, 
introduced by \cite{Young-SD-1917} and \cite{Daniell-SD-1918}, 
which under obvious restrictions on $g(t)$ (for details see, e.g., \cite{Kim-SD-2011}), 
can be written as
\begin{equation}\label{Stieltjes-derivative}
	\frac{df(t)}{dg(t)} = 
	\lim_{|t_1 - t_2| \rightarrow 0}   
	\frac{f(t_1) - f(t_2)}
	       {g(t_1) - g(t_2)}, 
\end{equation} 

\begin{figure}[h!]
\begin{center}
\includegraphics[width=8cm]{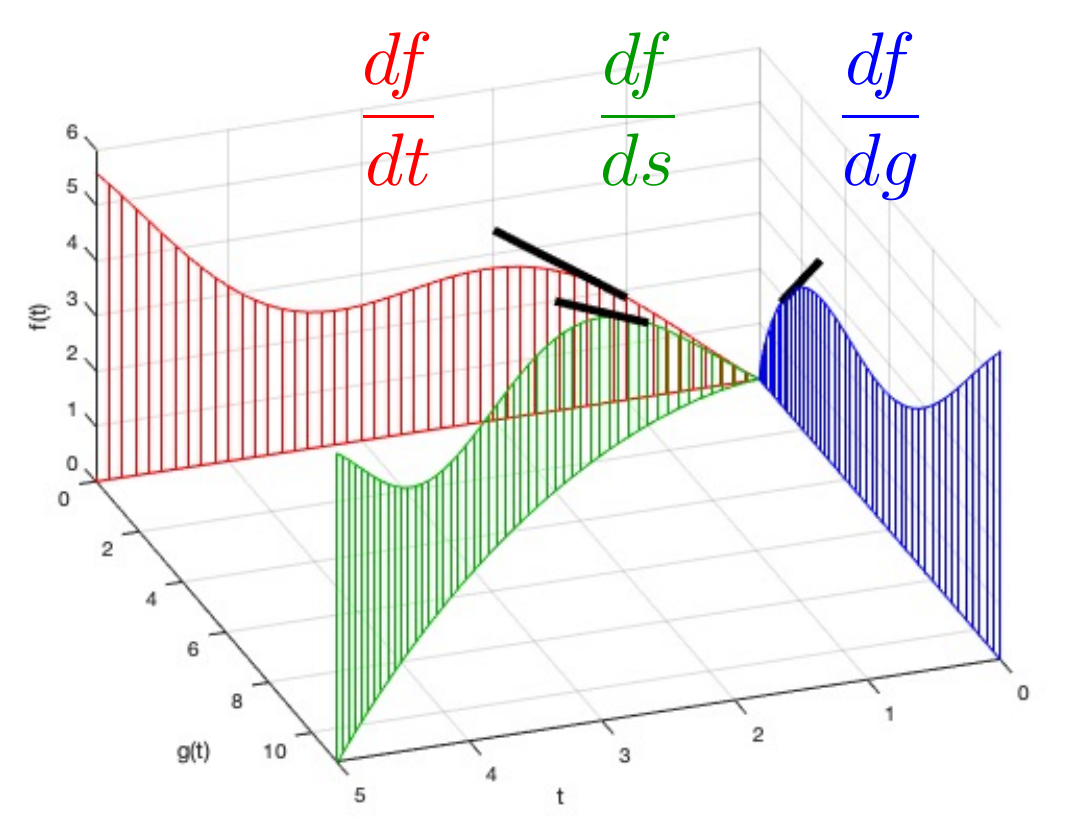}    
\caption{Interpretation of the classical integral (red area), the integral along the path (green area) and
Stieltjes integral (blue area). The tangent lines serve as interpretations of the classical first-order derivative (red), of the first-order derivative along the path (green), and of the Stieltjes derivative (blue).} 
\label{fig:Stieltjes-interpretation}
\end{center}
\end{figure}

Looking now at Fig.~\ref{fig:Stieltjes-interpretation}, we have the geometric meaning 
of the Stieltjes derivative (\ref{Stieltjes-derivative}): 
it is the slope of the tangent line at the point $( g(t), f(t) )$ in the ``blue'' plane $(\tau, y)$. 

The classical first-order derivative is the slope of the tangent line at the point $(t, f(t) )$ 
in the ``red'' plane $(t, y)$.

Finally, we can introduce the derivative along the path $\tau = g(t)$
(note that is is not a directional derivative) as
\begin{equation}\label{path-derivative}
		\frac{df(t)}{ds(t)} = 
	\lim_{|t_1 - t_2| \rightarrow 0}   
	\frac{f(t_1) - f(t_2)}
	       {s(t_1) - s(t_2)}, 
\end{equation}
where $s(t_k)$ is the length of the curve $\tau = g(t)$ for the interval $a \leq t \leq t_k $ ($k = 1, 2$). 
The geometric meaning of the derivative along the path (\ref{path-derivative}) is the slope 
of the tangent line to the ``green fence'' at the point $(t, g(t), f(t))$. 

\section{Geometric interpretation of the~``fractal derivative''}

The notion of the ``fractal derivative'' 

\vspace*{0.3ex}
\begin{equation}\label{path-derivative}
		\frac{df(t)}{dt^\alpha} = 
	\lim_{t_1 \rightarrow  t}   
	\frac{f(t_1) - f(t)}
	       {t_1^\alpha - t^\alpha},
	 \quad 
	 \alpha > 0, 
\end{equation}
\vspace*{-1ex}
that was first introduced by \cite{{Wen-Chen-2006}}, is obviously just a particular case 
of the Stieltjes derivative (\ref{Stieltjes-derivative}) for $g(t) = t^\alpha$. 
Thus the geometric meaning of the ``fractal derivative'' 
is the same as of Stieltjes derivative -- the slope of the tangent line 
at the point $(t^\alpha, f(t))$  in the ``blue'' plane $(\tau, y)$ 
shown in Fig.~\ref{fig:Stieltjes-interpretation}.


\vspace*{2cm}

\section*{Place and date of the talk}
This work has been submitted as discussion paper 074 and presented as talk WeB02.2 (Wednesday, July 9, 2024, 16:40-17:00) at \textit {ICFDA'2024  -- 12th IFAC Conference on Fractional Differentiation and its Applications, Bordeaux, France, July 9--12, 2024.}

https://icfda2024.sciencesconf.org/

https://ifac.papercept.net/conferences/scripts/rtf/FDA24\_ContentListWeb\_2.html

                                                   






\end{document}